\documentclass[10pt]{amsart}
\input{epsf}
\usepackage{amssymb,latexsym}
\topskip  \numberwithin{equation}{section}
\newtheorem{theorem}{Theorem}[section]
\newtheorem{proposition}[theorem]{Proposition}
\newtheorem{corollary}[theorem]{Corollary}

\newtheorem{example}[theorem]{Example}

\begin{document}

\title{ {\sc restricted $132$-involutions and Chebyshev polynomials}}

\author[Olivier Guibert and Toufik Mansour]
{O. Guibert and T. Mansour } \maketitle
\begin{center}
{\small LABRI (UMR 5800), Universit\'e Bordeaux 1,
       351 cours de la Lib\'eration, 33405 Talence Cedex, France}\\[4pt]
{\tt guibert@labri.fr}, {\tt toufik@labri.fr}
\end{center}

\section*{Abstract}

We study generating functions for the number of involutions in
$S_n$ avoiding (or containing once) $132$, and avoiding (or
containing once) an arbitrary permutation $\tau$ on $k$ letters.
In several interesting cases the generating function depends only
on $k$ and is expressed via Chebyshev polynomials of the second
kind. In particular, we establish that involutions avoiding both
$132$ and $12\dots k$ have the same enumerative formula according
to the length than involutions avoiding both $132$ and any {\em
double-wedge pattern} possibly followed by fixed points of total
length $k$. Many results are also shown with a combinatorial point
of view.

\section{Introduction}
A permutation is a bijection from $[n]=\{1,2,\ldots,n\}$ to $[n]$.
Let $S_n$ be the set of permutations of length $n$.

Let $\alpha\in S_n$ and $\tau\in S_k$ be two permutations. We say
that $\alpha$ {\it contains\/} $\tau$ if there exists a
subsequence $1\leq i_1<i_2<\dots<i_k\leq n$ such that
$(\alpha_{i_1},\dots,\alpha_{i_k})$ is order-isomorphic to $\tau$;
in such a context $\tau$ is usually called a {\it pattern\/}. We
say that $\alpha$ {\it avoids\/} $\tau$, or is $\tau$-{\it
avoiding\/}, if such a subsequence does not exist. The set of all
$\tau$-avoiding permutations in $S_n$ is denoted $S_n(\tau)$. For
an arbitrary finite collection of patterns $T$, we say that
$\alpha$ avoids $T$ if $\alpha$ avoids any $\tau\in T$; the
corresponding subset of $S_n$ is denoted $S_n(T)$.

While the case of permutations avoiding a single pattern has attracted
much attention, the case of multiple pattern avoidance remains less
investigated. In particular, it is natural, as the next step, to consider
permutations avoiding pairs of patterns $\tau_1$, $\tau_2$. This problem
was solved completely for $\tau_1,\tau_2\in S_3$ (see \cite{SS}), for
$\tau_1\in S_3$ and $\tau_2\in S_4$ (see \cite{W}), and for
$\tau_1,\tau_2\in S_4$ (see \cite{B1,Km} and references therein).
Several recent papers \cite{CW,MV1,Kr,MV2} deal with the case
$\tau_1\in S_3$, $\tau_2\in S_k$ for various pairs $\tau_1,\tau_2$. Another
natural question is to study permutations avoiding $\tau_1$ and containing
$\tau_2$ exactly $t$ times. Such a problem for certain $\tau_1,\tau_2\in S_3$
and $t=1$ was investigated in \cite{R}, and for certain $\tau_1\in S_3$,
$\tau_2\in S_k$ in \cite{RWZ,MV1,Kr}. The tools involved in these papers
include continued fractions, Chebyshev polynomials, and Dyck words.

{\it Chebyshev polynomials of the second kind\/} (in what follows
just Chebyshev polynomials) are defined by
    $$U_r(\cos\theta)=\frac{\sin(r+1)\theta}{\sin\theta}$$
for $r\geq0$. Evidently, $U_r(x)$ is a polynomial of degree $r$ in
$x$ with integer coefficients. Chebyshev polynomials were invented
for the needs of approximation theory, but are also widely used in
various other branches of mathematics, including algebra,
combinatorics, and number theory (see \cite{Ri}).

{\it Dyck words} are words $w$ of $\{x,\overline{x}\}^*$ verifying
that $|w|_x = |w|_{\overline{x}}$ and that for all $w = w' w''$,
$|w'|_x \geq |w'|_{\overline{x}}$. Dyck words of length $2n$ are
enumerated by the $n$th Catalan number
$C_n=\frac{1}{n+1}\binom{2n}{n}$ whose generating function is
$C(x)=\frac{1-\sqrt{1-4x}}{2x}$.

We also consider words of $\{a,b^2\}^*$ of length $n$
enumerated by the $n$th Fibonacci number $F_n$
with $F_0=F_1=1$ and $F_{n}=F_{n-1}+F_{n-2}$
whose generating function is $F(x)=\frac{1}{1-x-x^2}$.

Apparently, for the first time the relation between restricted
permutations and Chebyshev polynomials was discovered  by Chow and
West in \cite{CW}, and later by Mansour and Vainshtein
\cite{MV1,MV2,MV3,MV4}, Krattenthaler \cite{Kr}. These results
related to a rational function
    $$R_k(x)=\frac{2t U_{k-1}(t)}{U_k(t)}, \qquad t=\frac{1}{2\sqrt{x}}$$
for all $k\geq 1$.\\

An involution is a permutation such that its cycles are of length
$1$ or $2$ that is $\alpha\in S_n$ is an involution if and only if
$\alpha(\alpha_i)=i$ for all $i\in [n]$.

Some authors considered involutions with forbidden patterns.
\\
Regev in \cite{Regev} provided asymptotic formula for $12\cdots
k$-avoiding involutions of length $n$ and he also established that
$1234$-avoiding involutions of length $n$ are enumerated by
Motzkin numbers. Gessel \cite{Gessel} exhibited the enumeration of
such $12\cdots k$-avoiding involutions of length $n$. Moreover,
Gouyou--Beauchamps \cite{GouyouBYoung} obtained by an entirely
bijective proof very nice exact formulas for the number of
$12345$-avoiding and $123456$-avoiding involutions of length $n$.
\\
Gire \cite{GireThese} studied some permutations with forbidden
subsequences and established a one-to-one correspondence between
1-2 trees having $n$ edges and permutations of length $n$ avoiding
patterns $321$ and $231$, the latter being allowed in the case
where it is itself a subsequence of the pattern $3142$. Guibert
\cite{GuibertThese} also established bijections between 1-2 trees
having $n$ edges and another set of permutations with forbidden
subsequences and $3412$-avoiding involutions of length $n$ and
$4321$-avoiding involutions of length $n$ (and so with
$1234$-avoiding involutions of length $n$ by transposing the
corresponding Young tableaux obtained by applying the
Robinson-Schensted algorithm). He also shown \cite{GuibertThese} a
bijection between vexillary involutions (that is $2143$-avoiding
involutions) and $1243$-avoiding involutions. More recently,
Guibert, Pergola and Pinzani \cite{GPP} established a one-to-one
correspondence between 1-2 trees having $n$ edges and vexillary
involutions of length $n$. So all these sets are enumerated by the
$n$th Motzkin number
$\sum_{i=0}^{\lfloor\frac{n}{2}\rfloor}\binom{n}{2i}C_{i}$. It
remains a connected open problem: in \cite{GuibertThese}
conjectures that $1432$-avoiding involutions of length $n$ are
also enumerated by the $n$th Motzkin number.

In this paper we present a general approach to the study of
involutions in $S_n$ avoiding $132$ (or containing $132$ exactly
once), and avoiding (or containing exactly once) an arbitrary pattern $\tau\in
S_k$. As a consequence, we derive all the previously known results
for this kind of problems, as well as many new results.
Some results are also established by bijections as for an example
a bijection between $132$-avoiding involutions and primitive Dyck words.

The paper is organized as follows. The case of involutions
avoiding both $132$ and $\tau$ is treated in Section $2$. We
present an explicit expression in terms of Chebyshev polynomials
for several interesting cases. The case of involutions avoiding
$132$ and containing $\tau$ exactly once is treated in section
$3$. Here again we present an explicit expression in terms of
Chebyshev polynomials for several interesting cases. Finally, the
cases of involutions containing $132$ exactly once and either
avoiding or containing exactly once an arbitrary pattern $\tau$ is
treated in sections $4$ and $5$; respectively.

\section{Avoiding $132$ and another pattern}
\label{sec2}

Let $I_T(n)$ denote the number of involutions in $S_n(132)$ avoiding
$T$, and let $I_T(x)=\sum_{n\geq 0}I_T(n)x^n$ be the corresponding
generating function. The following proposition is the base of all the other
results in this section, which holds immediately from definitions.


\begin{proposition}\label{prom}
For any involution $\pi\in S_n(132)$ such that $\pi_j=n$ holds
either,
\begin{enumerate}
\item   for $1\leq j\leq [n/2]$, $\pi=(\beta,n,\gamma,\delta,j)$, where
    \begin{itemize}
    \item[$\ \ I.$] $\beta$ is a $132$-avoiding permutation of the numbers $n-j+1,\dots,n-2,n-1$,

    \item[$\ II.$] $\delta$ is a $132$-avoiding permutation of the numbers $1,\dots,j-2,j-1$ such
    that $\delta\cdot\beta$ is the identity permutation of $S_{j-1}$,

    \item[$III.$] $\gamma$ is a $132$-avoiding involution of the numbers $j+1,j+2,\dots,n-j-1,n-j$;
    \end{itemize}

\item   for $j=n$, $\pi=(\beta,n)$ where $\beta$ is an involution in $S_{n-1}(132)$.
\end{enumerate}
\end{proposition}

As a corollary of Proposition \ref{prom} we get the generating
function for the number of involutions in $S_n(132)$ as follows.

\begin{theorem} {\rm (see \cite[Prop.~5]{SS})}
Let $C(t)$ be the generating function for
the Catalan numbers; then
                  $$I_\varnothing(x)=\frac{1}{1-x-x^2C(x^2)}.$$
\end{theorem}
\begin{proof}
Proposition \ref{prom} yields for all $n\geq 1$,
   $$I_\varnothing(n)=\sum_{j=1}^{[n/2]} C_{j-1}I_\varnothing(n-2j)+I_\varnothing(n-1),$$
where $C_{j-1}$ is the $(j-1)$th Catalan number. Besides
$I_\emptyset(0)=1$, therefore in terms of generating function we
get that
   $$I_\varnothing(x)=1+x^2C(x^2)I_\varnothing(x)+xI_\varnothing(x).$$
\end{proof}

We can also prove this result by a bijective point of view.


Let $P_{x,\overline{x}} = \{ w \in \{x,\overline{x}\}^* :
\mbox{ for all } w = w' w'', |w'|_x \geq |w'|_{\overline{x}} \}$
be the language of primitive Dyck words.
\\
The number of such words of $P_{x,\overline{x}}$ of length $n$ is
the central binomial coefficient
$\binom{n}{\lfloor\frac{n}{2}\rfloor}$ with $n \geq 0$. Indeed,
any primitive Dyck word $w$ of $P_{x,\overline{x}}$ can be
uniquely written as $w_0 x w_1 x \ldots x w_p$ where $w_i$ is a
Dyck word (that is $w_i \in P_{x,\overline{x}}$ and $|w_i|_x =
|w_i|_{\overline{x}}$) for all $0 \leq i \leq p$, but $w$ can also
be uniquely written as $w_0 \overline{x} w_1 \overline{x} \dots
\overline{x} w_{\lceil \frac{p}{2} \rceil - 1} \overline{x}
w_{\lceil \frac{p}{2} \rceil} x w_{\lceil \frac{p}{2} \rceil + 1}
x \ldots x w_p$. So primitive Dyck words $w$ of
$P_{x,\overline{x}}$ of length $n$ are in bijection with bilateral
words of $\{ w \in \{x,\overline{x}\}^* : |w|_x =
|w|_{\overline{x}} \mbox{ or } |w|_x = |w|_{\overline{x}} - 1 \}$
of length $n$ trivially enumerated by
$\binom{n}{\lfloor\frac{n}{2}\rfloor}$.


\begin{theorem}
\label{Phi} There is a bijection $\Phi$ between involutions
in $S_n(132)$ and primitive Dyck words of
$P_{x,\overline{x}}$ of length $n$. Moreover, the number of fixed
points of the involution corresponds to the difference between the
number of letters $x$ and $\overline{x}$ into the primitive Dyck
word.
\end{theorem}
\begin{proof}
Let $\pi$ be a $132$-avoiding involution on $[n]$ having $p$ fixed
points. According to Proposition~\ref{prom} we have
$\pi=\pi'\pi''x\pi'''$ with $|\pi'|=\frac{n-p}{2}$ ($\pi'$ has no
fixed points and constitutes cycles with $\pi''$ or $\pi'''$),
$\pi''$ does not contain fixed point and $\pi(x)=x$ ($x$ is the
first fixed point). We obtain two $132$-avoiding involutions on
$[n+1]$ from $\pi$: the first one is given by inserting a fixed
point between $\pi'$ and $\pi''$, and the second one (iff $\pi$
has at least one fixed point) is given by modifying the first
fixed point $x$ by a cycle starting between $\pi'$ and $\pi''$.
All $132$-avoiding involutions can be obtained (and only once) by
applying this rule, starting from the empty involution of length
$0$.
\\
Let $w$ be a primitive Dyck word of $P_{x,\overline{x}}$ of length
$n$ such that $|w|_x-|w|_{\overline{x}}=p$. So we have $w = w_0 x
w_1 x \ldots x w_p$ where $w_i$ are Dyck words for all $0 \leq i
\leq p$. We obtain two primitive Dyck words of length $n+1$ from
$w$: $x w$ and $x w_0 \overline{x} w_1 x \cdots x w_p$ (iff
$p>0$). All primitive Dyck words can be obtained (and only once)
by applying this rule, starting from the empty word of length $0$.
\end{proof}

Clearly, these two generating trees for the $132$-avoiding involutions and
the primitive Dyck words can be characterized by the following succession
system:\\
    $\left\{\begin{array}{lcll}
    \multicolumn{4}{l}{(0)} \\
    (0) & \leadsto & (1) & \\
    (p) & \leadsto & (p+1) , (p-1) & \mbox{if } p \geq 1
    \end{array}\right.$

\begin{figure}
\vspace*{-7.00truecm} \hspace*{-2.75truecm}
\epsfxsize=590.0pt 
\epsffile{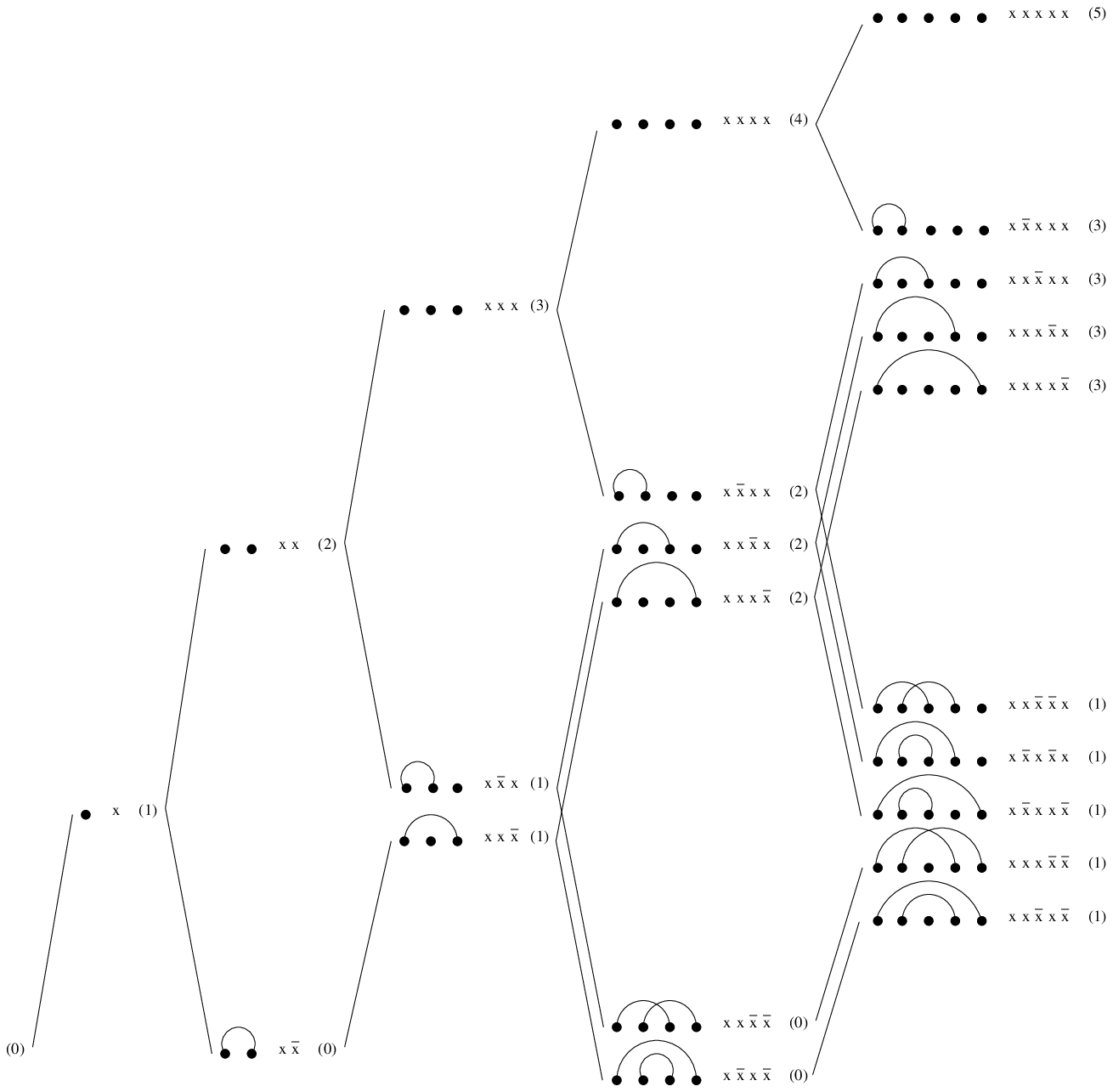} 
\vspace*{-7.75truecm} {\caption{The generating trees of the
$132$-avoiding involutions and the primitive Dyck words, with the
labels of the succession system which characterizes them.}
\label{fig-ag}}
\end{figure}

Figure~\ref{fig-ag} shows the bijection $\Phi$ between
$132$-avoiding involutions and primitive Dyck words
(and the labels of the succession system which characterizes them)
for the first values.

\begin{corollary}
\label{nb132ptfix} The number of $132$-avoiding involutions of
length $n$ is $\binom{n}{\lfloor \frac{n}{2} \rfloor}$. Moreover,
the number of $132$-avoiding involutions of length $n$ having $p$
fixed points with $0 \leq p \leq n$ (and $p$ is odd iff $n$ is
odd) is the ballot number $\binom{n}{\frac{n+p}{2}} -
\binom{n}{\frac{n+p}{2}+1}$.
\end{corollary}

\begin{proof}
Indeed, the number of primitive Dyck words $w$ of $P_{x,\overline{x}}$
according to the length and $|w|_x-|w|_{\overline{x}}$ is given by
the ballot number (or Delannoy number \cite{Errera} or
distribution $\alpha$ of the Catalan number \cite{Kre}).
\end{proof}

In particular, the number of $132$-avoiding involutions of length $2n$
without fixed points is $C_n$ the $n$th Catalan number.
%


The following theorem is the base of all the other results in this
section.

\begin{theorem}\label{thg}
Let $T$ set of patterns, $T'=\{(\tau,|\tau|+1) : \tau\in T\}$, and
let $S_T(x)$ be the generating function for the number of
$T$-avoiding permutations in $S_n(132)$. Then
    $$I_{T'}(x)=\frac{1}{1-x^2S_T(x^2)}+\frac{x}{1-x^2S_T(x^2)}I_T(x).$$
\end{theorem}
\begin{proof}
Proposition \ref{prom} with definitions of $T'$ yields for $n\geq
1$,
    $$I_{T'}(n)=I_T(n-1)+\sum_{j=1}^{[n/2]} s_T(j-1)I_{T'}(n-2j),$$
where $s_T(j-1)$ is the number of permutations in
$S_{j-1}(132,T)$. Hence, in terms of generating functions we have
    $$I_{T'}(x)-1=x\cdot I_T(x)+x^2 S_T(x^2)\cdot\frac{I_{T'}(x)+I_{T'}(-x)}{2}+x^2S_T(x^2)\cdot\frac{I_{T'}(x)-I_{T'}(-x)}{2},$$
so the theorem holds.
\end{proof}


\subsection{Avoiding $132$ and $12\dots k$}
\label{subsec21}

\begin{example}\label{ex1} {\rm (see \cite{SS})}
Let us find $I_{123}(x)$; let $T'=\{123\}$ and $T=\{12\}$,
so Theorem \ref{thg} gives
    $$I_{123}(x)=\frac{1}{1-x^2S_{12}(x^2)}+\frac{x}{1-x^2S_{12}(x^2)}I_{12}(x),$$
where by definitions $S_{12}(x)=I_{12}(x)=\frac{1}{1-x}$, hence
    $$I_{123}(x)=\frac{1+x}{1-2x^2},$$
which means the number of involutions $I_{123}(n)$ is given by
$2^{[n/2]}$ for all $n\geq 0$.

Similarly, $I_{1234}(x)=\frac{1}{1-x-x^2}$, so the number of
involutions $I_{1234}(n)$ is given by $F_n$ the $n$th Fibonacci number.
\end{example}

The case of varying $k$ is more interesting. As an extension of
Example~\ref{ex1} let us consider the case $T=\{12\dots k\}$.

\begin{theorem}\label{th12k}
For all $k\geq 1$,
    $$I_{12\dots k}(x)=\frac{1}{x\cdot U_k\left( \frac{1}{2x} \right)}
          \sum_{j=0}^{k-1} U_j\left( \frac{1}{2x} \right).$$
\end{theorem}
\begin{proof}
Immediately, the theorem holds for $k=1$. Let $k\geq 2$; Theorem
\ref{thg} gives
    $$I_{12\dots k}(x)=\frac{1}{1-x^2S_{12\dots (k-1)}(x^2)}+\frac{x}{1-x^2S_{12\dots (k-1)}(x^2)}I_{12\dots (k-1)}(x).$$
On the other hand, the generating function for the sequence
$S_n(132,12\dots(k-1))$ is given by $R_{k-1}(x)$ (see
\cite[Th.~1]{CW}) with $R_k(x)=\frac{1}{1-xR_{k-1}(x)}$ (see
\cite{MV1}) we get that
    $$I_{12\dots k}(x)=R_k(x^2)+xR_k(x^2)I_{12\dots (k-1)}(x).$$
Besides $I_{1}(x)=R_0(x)=1$, hence by use induction on $k$ and
definitions of $R_k(x)$ the theorem holds.
\end{proof}

We consider now a combinatorial point of view for this result.


Let $\pi$ be a $132$-avoiding involution. Clearly, if $\pi$ avoids
$12\dots k$ then $\pi$ has less than $k$ fixed points. Moreover,
if $\pi$ of length $n$ having less than $k$ fixed points is
obtained from an $132$-avoiding involution $\sigma$ of length less
than $n$ having $k$ fixed points (and take $\sigma$ as big as
possible) by applying the rules described for bijection $\Phi$
given by Theorem~\ref{Phi}, then $\pi$ contains a subsequence $1 2
\ldots k$ because the first fixed points of $\sigma$ become cycles
into $\pi$ such that the beginning of these cycles and the last
remaining fixed points of $\sigma$ into $\pi$ constitute a
subsequence of type $1 2 \ldots k$.
So the succession system $(*)$ \\
    $\left\{\begin{array}{lcll}
    (0) \\
    (0) & \leadsto & (1) \\
    (p) & \leadsto & (p+1) , (p-1) & 1 \leq p \leq k-2 \\
    (k-1) & \leadsto & (k-2)
    \end{array}\right.$\\
characterizes the generating tree of the involutions avoiding both
$132$ and $1 2 \ldots k$.

It is easy to see that for $k$ odd, the number of involutions of
length $2m$ avoiding both $132$ and $12\dots k$ is the twice of
the number of involutions of length $2m-1$ avoiding both $132$ and
$12\dots k$.

Moreover, the reader can note that the set of labels of this
succession system is finite and so the corresponding generating
function is rational. More precisely, we immediately deduce from
the previous succession system that the number of involutions of
length $n$ avoiding both $132$ and $12\dots k$ and having $p$
fixed points is given by the $(p+1)$th component of the vector
given by $V_k . M_k^n$ where $V_k = ( \ 1 \ 0 \ 0 \ \dots \ 0 \ )$
is a vector of $k$ elements and $M_k = \left(
\begin{array}{ccccccccc}
    0 & 1 & 0 & 0 & 0 & \cdots & 0 & 0 & 0 \\
    1 & 0 & 1 & 0 & 0 & \cdots & 0 & 0 & 0 \\
    0 & 1 & 0 & 1 & 0 & \cdots & 0 & 0 & 0 \\
    0 & 0 & 1 & 0 & 1 & \cdots & 0 & 0 & 0 \\
      &   &   &   &   & \ddots &   &   &   \\
    0 & 0 & 0 & 0 & 0 & \cdots & 1 & 0 & 1 \\
    0 & 0 & 0 & 0 & 0 & \cdots & 0 & 1 & 0 \\
\end{array} \right)$ is a $k \times k$ matrix.
\\
In another way, we can see that as an automaton where the states
are $0, 1, \dots, k-1$ and the transitions are arrows from $i$ to
$i+1$ for $0 \leq i < k-1$ and from $i$ to $i-1$ for $0<i\leq
k-1$.

The bijection $\Phi$ establishes a one-to-one correspondence
between involutions of length $n$ avoiding both $132$ and $12\dots
k$ and having $p$ fixed points, and primitive Dyck words $w = w_0
x w_1 x \ldots x w_p$ of $P_{x,\overline{x}}$ of length $n$ such
that $w_i$ is a Dyck word of height less than $k-p+i$ (that is
$w_i \in P_{x,\overline{x}}$ and $|w_i|_x = |w_i|_{\overline{x}}$
and for all $w_i = w' w''$, $|w'|_x-|w'|_{\overline{x}}<k-p+i$)
for all $0 \leq i \leq p$.
\\
In particular involutions of length $2n$ avoiding both $132$ and
$12\dots k$ without fixed points are in bijection by $\Phi$ with
Dyck words of length $2n$ of height less than $k$.

These Dyck words of bounded height was considered by Kreweras \cite{Kre}
and Viennot \cite{V}.
In particular, Dyck words of length $2n$ of
height less than $1$, $2$, $3$, $4$, $5$ are respectively enumerated by
$0$, $1$, $2^{n-1}$, $F_{n-2}$, $\frac{3^{n-1}+1}{2}$ for all $n \geq 1$.
\\
We provide some simple bijections for special cases $k=3,4,5$
(related to Example~\ref{ex1}) by generating some well known words
in the same way as involutions avoiding both $132$ and $1 2 \ldots k$.

Fist of all, we consider the case $k=3$ and the words of
$\{a,b\}^*$ or $a\{a,b\}^*$ enumerated by the powers of $2$
we can generate from the empty word labeled $(0)$ by the rules:\\
    $\left\{\begin{array}{lcl}
    w (0) & \leadsto & a w (1) \\
    a w (1) & \leadsto & a w (2) , b w (0) \\
    w (2) & \leadsto & a w (1)
    \end{array}\right.$\\
such that the words labeled $(0)$ start by $b$ whereas the words
labeled $(1)$ or $(2)$ start by $a$.
\\
So, words of $\{a,b\}^n$ (respectively $a\{a,b\}^n$)
are in bijection with involutions avoiding both $132$ and $123$
of length $2n$ (respectively $2n+1$)
enumerated by $2^n$ (respectively $2^n$).

Next we consider the case $k=4$ and the words of $\{a,b^2\}^\ast$
enumerated by the Fibonacci numbers we can generate
from the empty word labeled $(0)$ by the rules:\\
    $\left\{\begin{array}{lcl}
    w (0) & \leadsto & a w (1) \\
    a w (1) & \leadsto & a a w (2) , b^2 w (0) \\
    a w (2) & \leadsto & b^2 w (3) , a a w (1) \\
    w (3) & \leadsto & a w (2)
    \end{array}\right.$\\
such that the words labeled $(0)$ or $(3)$ start by $b^2$ whereas
the words labeled $(1)$ or $(2)$ start by $a$.
\\
So, words of $\{a,b^2\}^\ast$ of length $n$ are in bijection with
involutions in $S_n(132,1234)$ enumerated by $F_n$ the $n$th
Fibonacci number.

Now we consider the case $k=5$ and the words of $\{a,b,c\}^\ast a$
or $\{a,b,c\}^\ast a \cup b \{a,b,c\}^\ast a$ enumerated by the
powers of $3$ we can generate
from the empty word labeled $(0)$ by the rules:\\
    $\left\{\begin{array}{lcl}
    w (0) & \leadsto & a w (1) \\
    w (1) & \leadsto & b w (2) , w (0) \\
    w = b w' = b b^* c w (2) & \leadsto & w (3) , a w' (1) \\
    w = b w' = b b^* a w (2) & \leadsto & c w' (3) , w (1) \\
    w (3) & \leadsto & w (4) , b w (2) \\
    w (4) & \leadsto & c w (3)
    \end{array}\right.$\\
such that the words labeled $(0)$ or $(1)$ start by $b^\ast a$,
the words labeled $(3)$ or $(4)$ start by $b^\ast c$, and the
words labeled $(2)$ start by $b$ (and have one letter more than
words labeled $(0)$ or $(4)$ at the same level).
\\
So, words of $\{a,b,c\}^n a$
(respectively $\{a,b,c\}^n a \cup b \{a,b,c\}^n a$)
are in bijection with involutions avoiding both $132$ and $12345$
of length $2n+1$ (respectively $2n+2$)
enumerated by $3^n$ (respectively $2.3^n$).

Figure~\ref{fig-forbid} (an output of the software \verb/forbid/
\cite{GuibertThese}) shows the first values for the number of
involutions avoiding both $132$ and $12\dots k$ for $3 \leq k \leq
5$ according to the number of fixed points.

\begin{figure}[ht]
\begin{center}
Involutions $\pi \in S_n(132,123)$ according to $|\{\pi(x)=x\}|$ for $1 \leq n \leq 15$
\begin{verbatim}
         =1 =2 =2 =4 =4 =8 =8 =16 =16 =32 =32 =64 =64 =128 =128
      2:  -  1  0  2  0  4  0   8   0  16   0  32   0   64    0
      1:  1  0  2  0  4  0  8   0  16   0  32   0  64    0  128
      0:  0  1  0  2  0  4  0   8   0  16   0  32   0   64    0
         1: 2: 3: 4: 5: 6: 7:  8:  9: 10: 11: 12: 13:  14:  15: [n]
\end{verbatim}
Involutions $\pi \in S_n(132,1234)$ according to $|\{\pi(x)=x\}|$ for $1 \leq n \leq 15$
\begin{verbatim}
         =1 =2 =3 =5 =8 =13 =21 =34 =55 =89 =144 =233 =377 =610 =987
      3:  -  -  1  0  3   0   8   0  21   0   55    0  144    0  377
      2:  -  1  0  3  0   8   0  21   0  55    0  144    0  377    0
      1:  1  0  2  0  5   0  13   0  34   0   89    0  233    0  610
      0:  0  1  0  2  0   5   0  13   0  34    0   89    0  233    0
         1: 2: 3: 4: 5:  6:  7:  8:  9: 10:  11:  12:  13:  14:  15: [n]
\end{verbatim}
Involutions $\pi \in S_n(132,12345)$ according to $|\{\pi(x)=x\}|$ for $1 \leq n \leq 15$
\begin{verbatim}
         =1 =2 =3 =6 =9 =18 =27 =54 =81 =162 =243 =486 =729 =1458 =2187
      4:  -  -  -  1  0   4   0  13   0   40    0  121    0   364     0
      3:  -  -  1  0  4   0  13   0  40    0  121    0  364     0  1093
      2:  -  1  0  3  0   9   0  27   0   81    0  243    0   729     0
      1:  1  0  2  0  5   0  14   0  41    0  122    0  365     0  1094
      0:  0  1  0  2  0   5   0  14   0   41    0  122    0   365     0
         1: 2: 3: 4: 5:  6:  7:  8:  9:  10:  11:  12:  13:   14:   15: [n]
\end{verbatim}
{\caption{The number of involutions avoiding both $132$ and $1 2 \ldots k$
for $3 \leq k \leq 5$ according to the number of fixed points.}
\label{fig-forbid}}
\end{center}
\end{figure}


\subsection{Avoiding $132$ and $213\dots k$}

\begin{example}\label{ex2}
Let us find $I_{213}(x)$; let $T'=\{213\}$ and $T=\{21\}$,
so Theorem \ref{thg} gives
    $$I_{213}(x)=\frac{1}{1-x^2S_{21}(x^2)}+\frac{x}{1-x^2S_{21}(x^2)}I_{21}(x),$$
where by definitions $S_{21}(x)=I_{21}(x)=\frac{1}{1-x}$, hence
    $$I_{213}(x)=\frac{1+x}{1-2x^2},$$
which means the number of involutions $I_{213}(n)$ is given by
$2^{[n/2]}$ for all $n\geq 0$.

Similarly, $I_{2134}(x)=\frac{1}{1-x-x^2}$, so the number of
involutions $I_{2134}(n)$ is given by $F_n$ the $n$th Fibonacci
number.
\end{example}

We can easily prove by a combinatorial way that $I_{213}(n)$
is given by $2^{[n/2]}$.
\\
An involution $\pi$ avoiding both $132$ and $213$ of length $n$
can be either $(n+1-i)(n+2-i)\dots n \pi' 12\dots i$ with $i\geq1$
or $12\dots n$ such that $\pi'$ is also an involution avoiding
both $132$ and $213$ (of length $n-2i$, and if we subtract $i$ to
each element). We code this recursive decomposition of an
involution $\pi$ avoiding both $132$ and $213$ by a word of
nonnegative integers formed by the successive positive numbers $i$
and whose last nonnegative integer is (the smallest integer of)
the number of the fixed points in $\pi$ divided by $2$. This
coding is clearly bijective.
\\
For example, the involutions $\epsilon$, $1$, $12$, $21$, $123$,
$321$, $1234$, $4231$, $3412$ and $4321$ are respectively coded by
$0$, $0$, $1$, $10$, $1$, $10$, $2$, $11$, $20$ and $110$.
Moreover, the involution $21 \linebreak[0] \ 19 \linebreak[0] \ 20
\linebreak[0] \ 16 \linebreak[0] \ 17 \linebreak[0] \ 18
\linebreak[0] \ 15 \linebreak[0] \ 14 \linebreak[0] \ 9
\linebreak[0] \ 10 \linebreak[0] \ 11 \linebreak[0] \ 12
\linebreak[0] \ 13 \linebreak[0] \ 8 \linebreak[0] \ 7
\linebreak[0] \ 4 \linebreak[0] \ 5 \linebreak[0] \ 6
\linebreak[0] \ 2 \linebreak[0] \ 3 \linebreak[0] \ 1$ in
$S_{21}(132,213)$ is coded by $123112$.
\\
Thus, involutions avoiding both $132$ and $213$ of length $n$ are
coded by words $w = w_1 w_2 \ldots w_{l-1} w_l$ with $l \geq 1$,
$w_j \geq 1$ for all $1 \leq j < l$, $w_l \geq 0$ and
$\sum_{j=1}^{l} w_j = \lfloor \frac{n}{2} \rfloor$. Trivially,
words $w$ are in bijection with words $w_1 w_2 \ldots w_{l-1}
(w_l+1)$ which are compositions of $\lfloor \frac{n}{2} \rfloor +
1$ into $l$ positive parts enumerated by $2^{\lfloor \frac{n}{2}
\rfloor}$.

The case of varying $k$ is more interesting. As an extension of
Example \ref{ex2} let us consider the case $T=\{213\dots k\}$.
Similarly as Theorem \ref{th12k} we have

\begin{theorem}\label{th21k}
For all $k\geq 1$,
    $$I_{213\dots k}(x)=\frac{1}{x\cdot U_k\left( \frac{1}{2x} \right)}
          \sum_{j=0}^{k-1} U_j\left( \frac{1}{2x} \right).$$
\end{theorem}

Therefore, Theorem~\ref{th12k} and Theorem~\ref{th21k} yields
$I_{123\dots k}(n)=I_{213\dots k}(n)$.
We establish a bijection for this result.

\begin{theorem}
\label{12k=2134k}
There is a bijection between
involutions avoiding both $132$ and $12\ldots k$ of length $n$ and
involutions avoiding both $132$ and $2134\ldots k$ of length $n$,
for any $k \geq 3$.
\\
Moreover, two involutions in bijection have
the same number of fixed points $p$ for all $0 \leq p \leq k-3$
whereas the involutions avoiding both $132$ and $12\ldots k$
having $k-2$ or $k-1$ fixed points correspond to
the involutions avoiding both $132$ and $2134\ldots k$
having $k-2$ or more fixed points.
\end{theorem}
\begin{proof}
In order to establish this result we consider a generating tree
for the involutions avoiding both $132$ and $2134\ldots k$
which is characterized by the same succession system $(*)$
given in Subsection~\ref{subsec21}
characterizing a generating tree for the involutions
avoiding both $132$ and $12\ldots k$.
\\
So let $\pi$ be an involution avoiding both $132$ and $2134\ldots k$
of length $n$ and let $q=|\{\pi(x)=x\}|$ be the number of fixed points
of $\pi$.
The label $(p)$ of $\pi$ is defined by $p=q$ if $q \leq k-3$
or by $p=k-2$ if $q \geq k-2$ and $(n+k) \bmod 2 = 0$
or by $p=k-1$ if $q \geq k-2$ and $(n+k) \bmod 2 = 1$.
Of course, the empty involution of length $0$ has label $(0)$.
We obtain $\sigma$ an involution avoiding both $132$ and $2134\ldots k$
of length $n+1$ by applying the following rules:
\begin{itemize}
\item
If $p \in [0,k-3]$, we have $\pi=\pi'\pi''$ with $|\pi'|=\frac{n-p}{2}$,
and then $\sigma$ obtained by inserting a fixed point between
$\pi'$ and $\pi''$ has label $(p+1)$.
\item
If $p=k-2$, we have $\pi=\pi'x\pi''$ with $\pi(x)=x=\frac{n+4-k}{2}$,
and then $\sigma$ obtained by inserting a fixed point between
$\pi'$ and $x$ has label $(k-1)$.
\item
If $p \in [1,k-3]$, we have $\pi=\pi'\pi''x\pi'''$ with
$|\pi'|=\frac{n-p}{2}$, $\pi(x)=x$ and
$\pi(y) \neq y$ for all $1 \leq y < x$.
Then $\sigma$ obtained by modifying the first fixed point $x$ by
a cycle starting between $\pi'$ and $\pi''$ (and ending in $x$)
has label $(p-1)$.
\item
If $p=k-1$, we have $\pi=\pi'x(x+1)\pi''$ with $\pi(x)=x=\frac{n+3-k}{2}$,
and then $\sigma$ obtained by inserting a fixed point between
$\pi'$ and $x$ has label $(k-2)$.
\item
If $p=k-2$, we have $\pi=\pi'(x-j)(x-j+1)\dots(x+j)\pi''$ with $j
\geq 0$, $\pi(x)=x=\frac{n+4-k}{2}$, $|\pi'|=\frac{n-k}{2}+1-j$
and $e>x$ for all $e \in \pi'$. Then $\sigma$ obtained by
modifying the $2j+1$ fixed points between $\pi'$ and $\pi''$ by
$j+1$ consecutive cycles each of difference (between the index and
the value) $j+1$ that is
$(\pi_1'+1)(\pi_2'+1)\dots(\pi_{\frac{n-k}{2}+1-j}'+1)
(\frac{n-k}{2}+3)(\frac{n-k}{2}+4)\dots(\frac{n-k}{2}+3+j)
(\frac{n-k}{2}-j+2)(\frac{n-k}{2}-j+3)\dots(\frac{n-k}{2}+2)
\pi''$ has label $(k-3)$.
\end{itemize}
\end{proof}

\begin{corollary}
There is a bijection between
permutations avoiding both $132$ and $12\ldots k$ of length $n$ and
permutations avoiding both $132$ and $2134\ldots k$ of length $n$,
for any $k \geq 3$.
\end{corollary}
\begin{proof}
By Proposition~\ref{prom}, we deduce that permutations $\pi$
avoiding both $132$ and $12\ldots k$ (respectively $2134\ldots k$)
of length $n$ are in bijection with involutions without fixed points
$(\pi^{-1}+n) \pi$ avoiding both $132$ and $12\ldots k$
(respectively $2134\ldots k$) of length $2n$.
Moreover, a particular case of Theorem~\ref{12k=2134k} establishes
a one-to-one correspondence between
involutions avoiding both $132$ and $12\ldots k$ without fixed points and
involutions avoiding both $132$ and $2134\ldots k$ without fixed points.
\end{proof}



\subsection{Avoiding $132$ and $(d+1(d+2)\dots k12\dots d$}

\begin{example}\label{ex22}
By Proposition \ref{prom} it is easy to obtain for $n\geq1$,
    $$I_{231}=n;\quad  I_{321}=[n/2]+1.$$
\end{example}


We consider a combinatorial approach to show Example~\ref{ex22}.
Clearly, we have that involutions avoiding both $132$ and $231$ of length $n$
are $i (i-1) \ldots 1 (i+1) (i+2) \ldots n$ for all $1 \leq i \leq n$
and that involutions avoiding both $132$ and $321$ of length $n$
are $(i+1) (i+2) \ldots (2i) 1 2 \ldots i (2i+1) (2i+2) \ldots n$
for all $0 \leq i \leq \lfloor \frac{n}{2} \rfloor$.

%

As an extension of Example~\ref{ex22} let us consider the case
$T=\{[k,d]\}$, where $[k,d]=(d+1,d+2,\dots,k,1,2,\dots,d)$.

\begin{theorem}\label{thkd}
For any $k\geq 2$, $k/2\geq d\geq 1$,
    $$I_{[k,d]}=\frac{1}{x(U_d(t)-U_{d-1}(t))}\left[
                U_{d-1}(t)+\frac{U_{k-2d-1}(t)}{U_{k-d}(t)U_{k-d-1}(t)}\sum_{j=0}^{k-d-1} U_j(t)
                \right],\quad t=\frac{1}{2x}.$$
\end{theorem}
\begin{proof}
Proposition \ref{prom} yields, in the second case the generating
function for the number of involutions $[k,d]$-avoiding
permutations is $xI_{[k,d]}(x)$. In the first case, we assume that
$\gamma$ either $(1)$ avoiding $12\dots (k-d)$, or $(2)$
containing $12\dots (k-d)$. In $(1)$, $\beta$ and $\delta$
avoiding $12\dots (k-d-1)$, so the generating function for these
number of involutions is $x^2R_{k-d-1}(x^2)I_{12\dots(k-d)}(x)$
(similarly Theorem \ref{th12k}). In $(2)$, $\beta$ and $\delta$
avoiding $12\dots (d-1)$, so the generating function for these
number of involutions is
$x^2R_{d-1}(x^2)(I_{[k,d]}(x)-I_{12\dots(k-d)}(x))$ (the
generating function for the number of involutions in
$S_n(132,[k,d])$ such containing $12\dots(k-d)$ is given
$I_{[k,d]}(x)-I_{12\dots(k-d)}(x)$). Therefore
$$\begin{array}{ll}
I_{[k,d]}(x)&=1+xI_{[k,d]}(x)+x^2R_{k-d-1}(x^2)I_{12\dots(k-d)}(x)+\\
        &+x^2R_{d-1}(x^2)(I_{[k,d]}(x)-I_{12\dots(k-d)}(x)),
\end{array}$$
which means that
$$I_{[k,d]}(x)=\frac{1}{1-x-x^2R_{k-d-1}(x^2)}\cdot(1+x^2I_{12\dots (k-d)}(x)(R_{k-d-1}(x^2)-R_{d-1}(x^2))).$$
Hence, by use the identities $R_k(x)=\frac{1}{1-xR_{k-1}(x)}$ and
$R_a(x)-R_b(x)=\frac{U_{a-b-1}(t)}{\sqrt{x}U_a(t)U_b(t)}$, the
theorem holds.
\end{proof}

\begin{example}
\label{ex3412}
Theorem \ref{thkd} yields for $k=4$ and $d=2$, the number of
involutions $I_{3412}(n)$ is given by $F_n$ the $n$th Fibonacci number.
\end{example}

We consider a combinatorial approach to show Example \ref{ex3412}.
An involution $\pi$ avoiding both $132$ and $3412$ of length $n$
can be written $i\pi' 1(i+1)(i+2)\cdots n$ with $1\leq i\leq n$
such that $\pi'$ is also an involution avoiding both $132$ and
$3412$ (of length $i-2$, and if we subtract $1$ to each element).
We code $\pi$ by a word of $\{a, b^2\}^\ast$ of length $n$ in that
way: $a$ if $\pi_i=i$, $b^2$ if $\pi_i<i$ and nothing if $\pi_i>i$
for all $1\leq i\leq n$. This coding is clearly bijective.\\

Following \cite{MV2} we say that $\tau\in S_k$ is a {\it wedge\/}
pattern if it can be represented as
$\tau=(\tau^1,\rho^1,\dots,\tau^r,\rho^r)$ so that each of
$\tau^i$ is nonempty, $(\rho^1,\rho^2,\dots,\rho^r)$ is a layered
permutation of $1,\dots,s$ for some $s$, and
$(\tau^1,\tau^2,\dots,\tau^r)=(s+1,s+2,\dots,k)$. For example,
$645783912$ and $456378129$ are wedge patterns.

For a further generalization of Theorem \ref{th12k}, Theorem
\ref{th21k} and \cite[Th.~2.6]{MV2}, consider the following
definition. We say that $\tau\in S_{2l}$ is a {\it double-wedge\/}
pattern if there exist a wedge pattern $\sigma\in S_{l-1}$ such
that
  $$\tau=(\sigma^{-1}+l,2l,\sigma,l)\ \mbox{or}\ \tau=(\sigma+l,2l,\sigma^{-1},l).$$
For example, the double-wedge patterns of length $10$ are
$6789(10)1234\linebreak[0] 5$, $7689(10)\linebreak[0]
213\linebreak[0] 4\linebreak[0] 5\linebreak[0]$,
$7869(10)\linebreak[0] 31245\linebreak[0]$, $7896(10)\linebreak[0]
41235\linebreak[0]$, $8679(10)\linebreak[0] 23145\linebreak[0]$,
$8796(10)\linebreak[0] 42135\linebreak[0]$, $8967(10)\linebreak[0]
34125\linebreak[0]$, $9678(10)\linebreak[0] 23415\linebreak[0]$
and $9768(10)\linebreak[0] 32415\linebreak[0]$.

\begin{theorem}\label{dwedge}
For any double-wedge pattern $\tau\in S_{2l}(132)$
  $$I_\tau(x)=I_{12\dots(2l)}(x)=\frac{R_l(x^2)}{1-xR_l(x^2)}.$$
\end{theorem}
\begin{proof}
First of all, let us find the generating function $I_\rho(x)$
where $\rho=(\sigma^{-1},2l,\sigma,l)$. By use Proposition
\ref{prom} we obtain in the first case $xI_\rho(x)$, and in the
second case $x^2S_\sigma(x^2)I_\rho(x)$ where $S_\sigma(x^2)$ is
the generating function for the number of permutations in
$S_n(132,\sigma)$, therefore ($1$ for the empty permutation)
$$I_\rho(x)=1+xI_\rho(x)+x^2S_\sigma(x^2)I_\rho(x).$$
On the other hand, Mansour and Vainshtein proved
$S_\sigma(x)=R_{l-1}(x)$ for any wedge pattern $\sigma$, so
$$I_\rho(x)=\frac{1}{1-x-x^2R_{l-1}(x^2)}.$$
By use the identity $R_l(x)=\frac{1}{1-xR_{l-1}(x)}$ we have
     $$I_\rho(x)=\frac{R_l(x)}{1-xR_l(x^2)}.$$
Now, let us find $I_{12\dots 2l}(x)$ in terms of $R_j(x)$. So, by
use the identity
$$\sum_{j=0}^{2l}U_j(t)=\frac{U_{2l}(t)U_{l-1}(t)}{U_l(t)-U_{l-1}(t)}$$
and use the symmetric inverse operation, the first part of the
theorem holds.
\end{proof}

\begin{theorem}\label{extwedge}
For any wedge pattern $\sigma\in S_{l-1}$ the generating function
for the number of permutations in
$S_n(132,(\sigma^{-1}+l,2l,\sigma,l,2l+1,\dots,k))$ (or
$S_n(132,(\sigma+l,2l,\sigma^{-1},l,2l+1,\dots,k))$, or
$S_n(132,(\sigma+l,2l,\sigma,l,2l+1,\dots,k))$) is given by
$R_{k}(x)$, for all $k\geq 2l$.
\end{theorem}
\begin{proof}
Let $\tau=(\sigma^{-1}+l,2l,\sigma,l)$ and let $S_\tau(x)$ be the
generating function for the number of permutations in
$S_n(132,\tau)$. By \cite[Th.~1]{MV2} we have
$$S_\tau(x)=1+x(S_\tau(x)-S_\sigma^{-1}(x))S_\sigma(x)+xS_\sigma^{-1}(x)S_\tau(x).$$
On the other hand, by \cite[Th.~2.6]{MV2} and $\sigma$ a wedge
pattern in $S_{l-1}(132)$ we have
$S_{\sigma^{-1}}(x)=S_\sigma(x)=R_{l-1}(x)$. Therefore, by use the
identity $R_l(x)=\frac{1}{1-xR_{l-1}(x)}$ we get
$$S_\tau(x)=\frac{R_l(x)(1-xR_{l-1}(x)R_l(x))}{1-xR_l^2(x)}.$$
By use the definitions of Chebyshev polynomials of the second kind
it is easy to see
    $$\frac{R_l(x)(1-xR_{l-1}(x)R_l(x))}{1-xR_l^2(x)}=R_{2l}(x),$$
hence by use again Theorem \cite[Th.~1]{MV2} we have
$S_{(\tau,2l+1,\dots,k)}(x)=R_k(x)$. Similarly we obtain the other
cases.
\end{proof}

As a corollary of Theorem \ref{dwedge} we have

\begin{theorem}\label{idwedge}
For any double wedge pattern $\tau\in I_{2l}(132)$
    $$I_{(\tau,2l+1,2l+2,\dots,k)}(x)=I_{12\dots k}(x).$$
\end{theorem}
\begin{proof}
Since, if $S_\beta(x)=S_\gamma(x)$ and $I_\beta(x)=I_\gamma(x)$,
then Theorem \ref{thg} yields $I_{\tau'}(x)=I_{\beta'}(x)$, and by
use \cite[Th. 1]{MV2} we have $S_{\tau'}(x)=S_{\rho'}(x)$, where
$\tau'=(\tau_1,\dots,\tau_p,p+1)$ and
$\rho'=(\rho_1,\dots,\rho_p,p+1)$ two patterns in $S_{p+1}$.
Hence, the theorem holds by use Theorem \ref{dwedge}, Theorem
\ref{extwedge}, and induction on $p$.
\end{proof}

%

In view of Theorem \ref{dwedge} and Theorem \ref{idwedge} it is a
challenge to find a bijective proof.


\subsection{Avoiding $132$ and two other patterns}

Now, let us restrict more than two patterns ($132$ and two other
patterns).

\begin{example}\label{ex3}
Let us find $I_{123,213}(x)$; let $T'=\{123,213\}$ and $T=\{12,21\}$,
so Theorem \ref{thg} gives
    $$I_{123,213}(x)=\frac{1}{1-x^2S_{12,21}(x^2)}+\frac{x}{1-x^2S_{12,21}(x^2)}I_{12,21}(x),$$
where by definitions $S_{12,21}(x)=I_{12,21}(x)=1+x$, hence
    $$I_{123,213}(x)=\frac{1+x+x^2}{1-x^2-x^4},$$
which means the number of involutions $I_{123,213}(2n)$ is given
by $F_{n+1}$, and $I_{123,213}(2n+1)$ is given by $F_{n}$ for
all $n\geq 0$, where $F_m$ is the $m$th Fibonacci number.
\end{example}
We consider a combinatorial approach to show Example~\ref{ex3} by
distinguishing the cases of odd and even length.
\\
An involution $\pi$ avoiding $132$, $123$ and $213$ of length
$2n+1$ can be written either $(2n+1) \pi' 1$ or $(2n) (2n+1) \pi''
2 1$ or $1$ (if $n=0$) such that $\pi'$ and $\pi''$ are also
involutions avoiding $132$, $123$ and $213$ (of length $2n-1$ for
$\pi'$ if we subtract $1$ to each element, of length $2n-3$ for
$\pi''$ if we subtract $2$ to each element). We code $\pi$ by a
word of $\{a,b^2\}^*$ of length $n$ in that way: $a$ if
$\pi_i=2n+2-i$, $b^2$ if $\pi_i=2n+1-i$ and nothing if
$\pi_i=2n+3-i$ for all $1\leq i\leq n$. This coding is clearly
bijective.
\\
An involution $\pi$ avoiding $132$, $123$ and $213$ of length $2n$
can be written either $(2n) \pi' 1$ (that includes $2 1$ if $n=1$)
or $(2n-1) (2n) \pi'' 2 1$ or $1 2$ (if $n=1$) or the empty
involution (if $n=0$) such that $\pi'$ and $\pi''$ are also
involutions avoiding $132$, $123$ and $213$ (of length $2n-2$ for
$\pi'$ if we subtract $1$ to each element, of length $2n-4$ for
$\pi''$ if we subtract $2$ to each element). We code $\pi$ by a
word of $\{a,b^2\}^*$ of length $n+1$ in that way: $a$ if
$\pi_i=2n+1-i$ for all $1\leq i\leq n-1$, $b^2$ if $\pi_n=n+1$,
$b^2$ if $\pi_i=2n-i$ for all $1\leq i\leq n-2$, $b^2a$ if
$\pi_{n-1}=n+1$ and $aa$ if $\pi_n=n$. Moreover, the empty
involution is coded by $a$. This coding is clearly bijective.

%
%
%

Using definitions and Theorem \ref{thg} it is easy to see the
following.

\begin{corollary}\label{th12k213}
For all $k\geq 1$,
    $$I_{123\dots k, 213}(x)=I_{(k-1)\dots 21k, 123}(x)=\frac{1+x+x^2+\dots+x^{k-1}}{1-x^2-x^4-\dots-x^{2(k-1)}}.$$
\end{corollary}

\begin{example}
\label{ex312decr} Using Proposition \ref{prom} it is easy to see
for $n\geq 1$,
$$I_{213,321}(n)=\frac{1}{2}((-1)^n+3),\quad
I_{213,4321}(n)=[n/2]+1.$$
\end{example}

We consider a combinatorial approach to show Example~\ref{ex312decr}.
Clearly, we have that involutions avoiding $132$, $213$ and $321$
of length $n$ are $1 2 \ldots n$ and also
$(m+1) (m+2) \ldots n 1 2 \ldots m$ if $n=2m$ with $m \geq 1$.
We also have that involutions avoiding $132$, $213$ and $4321$ of length $n$ are
$(n+1-i) (n+2-i) \ldots n (i+1) (i+2) \ldots (n-i) 1 2 \ldots i$
for all $0 \leq i \leq \lfloor \frac{n}{2} \rfloor$.

%

\section{Avoiding $132$ and containing another pattern}
Let $I_\tau^r(n)$ denote the number of involutions in $S_n(132)$
containing $\tau$ exactly $r$ times, and let
$I_\tau^r(x)=\sum_{n\geq 0}I_\tau^r(n)x^n$ be the corresponding
generating function. Let us start by the following example.

\begin{example}\label{exbb0}
By Proposition \ref{prom} it is easy to see
$$I_{12}^1(x)=xI_{1}^1(x)+x^2I_{12}^1(x),$$
which means $I_{12}^1(x)=\frac{x^2}{1-x^2}$.
\end{example}

As extension of Example \ref{exbb0} let us consider the case
$\tau=12\dots k$.

\begin{theorem}\label{th1_12}
For all $k\geq 1$;
    $$I_{12\dots k}^{1}=\frac{1}{U_k\left( \frac{1}{2x} \right)}.$$
\end{theorem}
\begin{proof}
By Proposition \ref{prom} we have for $n\geq k$,
    $$I_{12\dots k}^{1}(n)=I_{12\dots (k-1)}^{1}(n-1)+\sum_{j=1}^{[n/2]} s_{12\dots (k-1)}(j-1)I_{12\dots k}^{1}(n-2j),$$
where $s_{12\dots k}(j-1)$ is the number of $12\dots k$-avoiding
permutations in $S_{j-1}(132)$. Besides $I_{12\dots k}^{1}(n)=0$
for all $n\leq k-1$, and $I_{12\dots k}^{1}(k)=1$. Similarly as
proof of Theorem \ref{th12k} we have
    $$I_{12\dots k}^{1}(x)=xR_{k}(x^2)I_{12\dots (k-1)}^{1}(x).$$
Hence, by induction on $k$ with initial condition $I_1^{1}=x$, the
theorem holds.
\end{proof}

Similarly as Theorem \ref{th1_12} we have an explicit formula when
$\tau=213\dots k$ or $\tau=23\dots k1$.

\begin{theorem}
For all $k\geq 2$;
    $$I_{213\dots k}^{1}=\frac{1-x^2}{U_k\left( \frac{1}{2x} \right)},
    \quad
    I_{23\dots k1}^1(x)=\frac{x^3}{(1-x)U_{k-2}\left( \frac{1}{2x} \right)}.$$
\end{theorem}

More generally, by Proposition \ref{prom} and the argument proof
of Theorem \ref{th12k} we get

\begin{theorem}\label{th12kr}
For any $k,r\geq 1$
    $$I_{12\dots k}^{r}(x)=xI_{12\dots (k-1)}^{r}(x)+x^2\sum_{2a+b=r} S_{12\dots (k-1)}^{a}(x^2)I_{12\dots k}^{b}(x),$$
where $S_{12\dots (k-1)}^{a}(x)$ is the generating function for
the number of permutations in $S_n$ containing $12\dots (k-1)$
exactly $a$ times.
\end{theorem}

In \cite{Kr} found an explicit formula for $S_{12\dots k}^{r}(x)$,
so Theorem \ref{th12kr} yields a recurrence for $I_{12\dots
k}^{r}(x)$. For example, following \cite{Kr} (\cite[Th. 3.1]{MV1})
we have a recurrence for $I_{12\dots k}^{r}(x)$ where
$r=1,2,\dots,2k$.

\begin{theorem}
Let $k\geq 1$; for $r=1,2,\dots,2k$
    $$I_{12\dots k}^{r}(x)=xI_{12\dots (k-1)}^r(x)+x^2R_{k-1}(x^2)I_{12\dots k}^{r}(x)
    +x^2\sum_{2a+b=r,\ a>0} x^{a-1}I_{12\dots k}^{b}(x)\frac{U_{k-1}^{a-1}\left( \frac{1}{2x} \right)}{U_k^{a+1}\left( \frac{1}{2x} \right)}.$$
\end{theorem}

The above Theorem yields for $r=2$ an explicit formula for
$I_{12\dots k}^{2}(x)$.

\begin{corollary}
For all $k\geq 1$,
    $$I_{12\dots k}^{2}(x)=\frac{1}{U_k\left( \frac{1}{2x} \right)}\sum_{i=1}^k
    \frac{\sum_{j=0}^{k-i} U_j\left( \frac{1}{2x} \right)}{U_{k+1-i}\left( \frac{1}{2x} \right) U_{k-i}\left( \frac{1}{2x} \right)}.$$
\end{corollary}

%
%

\section{Containing $132$ once and avoiding another pattern}

We first relate
involutions containing $132$ once to $132$-avoiding involutions.

\begin{theorem}
\label{bij132once-132} There is a bijection $\Psi$ between
involutions containing $132$ exactly once of length $n$ having $p$
fixed points with $1\leq p\leq n$ and $132$-avoiding involutions
of length $n-2$ having also $p$ fixed points.
\end{theorem}

\begin{proof}
Let $\pi=\pi'xz\pi''y\pi'''$ with $\pi(x)=x$, $\pi(y)=z$ and
$1+x=y<z$ be an involution containing $132$ once (that is
subsequence $xzy$) of length $n$ having $p$ fixed points. We
replace the subsequence $xzy$ by a fixed point between $\pi''$ and
$\pi'''$ in order to obtain an $132$-avoiding involution of length
$n-2$ having $p$ fixed points. Note that the only possibility to
have exactly once $132$ subsequence is a cycle with a fixed point
just to its left. Moreover, $y=x+1$ in order to forbid another
$132$ subsequence and cycles are only allowed from $\pi'$ to
$\pi''$ and from $\pi'$ to $\pi'''$ (and not from $\pi''$, $\pi'$,
$\pi''$, $\pi'''$ respectively to $\pi'''$, $\pi'$, $\pi''$,
$\pi'''$) whereas fixed points can uniquely be into $\pi'''$.
Clearly the involution we obtain avoids $132$ and in particular,
the fixed point $z-2$ cannot be a part of an $132$-subsequence
because it cannot be the $3$ or $2$ (all the elements on its left
are greater than it) and it cannot be the $1$ (there is no cycle
starting on its right).
\\
Let $\sigma = \sigma' \sigma'' \sigma''' t \sigma''''$ with
$\sigma(t) = t$ and $\sigma(i) \neq i$ for all $1 \leq i < t$
(that is $t$ is the first fixed point), $\sigma'(i) > t$ for all
$1 \leq i \leq |\sigma'|$ (all the elements of $\sigma'$ are
cycles ending into $\sigma''''$), $\sigma''(i) \in
[|\sigma'\sigma''|+1,t-1]$ for all $1 \leq i \leq |\sigma''|$ and
$\sigma'''(i) \in [|\sigma'|+1,|\sigma'\sigma''|]$ for all $1 \leq
i \leq |\sigma'''|$ ($\sigma''\sigma'''$ is entirely constituted
by cycles from $\sigma''$ to $\sigma'''$) be an $132$-avoiding
involution of length $n-2$ having $p$ fixed points. We modify the
fixed point $t$ by a cycle starting between $\sigma''$ and
$\sigma'''$ (and ending between $\sigma'''$ and $\sigma''''$) and
by adding a fixed point just to the right of $\sigma''$ in order
to obtain an involution containing $132$ once of length $n$ having
$p$ fixed points. Proposition~\ref{prom} leads immediately to the
decomposition of $\sigma$. The involution we obtain contains $132$
exactly once that is the subsequence we modify and insert. There
is no other $132$-subsequence and in particular, the fixed point
inserted and the start of the new cycle cannot be the $3$ or $2$
of another $132$-subsequence (all the elements on their left are
greater than them), the fixed point inserted and the start of the
new cycle and the end of the new cycle cannot be the $1$ of
another $132$-subsequence (there is no cycle starting on their
right), the end of the new cycle cannot be the $3$ of another
$132$-subsequence (because in that case the $2$ must be connected
to $\sigma'$ and the $1$ must be the fixed point inserted or an
element of $\sigma'''$ that forms an $231$-subsequence), and the
end of the new cycle cannot be the $2$ of another
$132$-subsequence (because in that case the $1$ must be an element
of $\sigma'\sigma''$ or the start of the new cycle and the $3$
must be the fixed point inserted or an element of $\sigma'''$ that
forms an $312$-subsequence excepted for the fixed point inserted
and the new cycle).
\\
So we have established a bijection
between $\pi$ an involution containing $132$ once
and $\sigma$ an $132$-avoiding involution
where $t=z-2$, $\pi'$ corresponds to $\sigma'\sigma''$,
$\pi''=\sigma'''$ and $\pi'''$ corresponds to $\sigma''''$.
\end{proof}

For example, the involution $22 \linebreak[0] \ 19 \linebreak[0] \
17 \linebreak[0] \ 18 \linebreak[0] \ 16 \linebreak[0] \ 12
\linebreak[0] \ 11 \linebreak[0] \ 13 \linebreak[0] \ {\bf 9}
\linebreak[0] \ {\bf 14} \linebreak[0] \ 7 \linebreak[0] \ 6
\linebreak[0] \ 8 \linebreak[0] \ {\bf 10} \linebreak[0] \ 15
\linebreak[0] \ 5 \linebreak[0] \ 3 \linebreak[0] \ 4
\linebreak[0] \ 2 \linebreak[0] \ 20 \linebreak[0] \ 21
\linebreak[0] \ 1 \linebreak[0] \ 23$ containing $132$ once (the
subsequence $9$~$14$~$10$) corresponds to the $132$-avoiding
involution $20 \linebreak[0] \ 17 \linebreak[0] \ 15 \linebreak[0]
\ 16 \linebreak[0] \ 14 \linebreak[0] \ 10 \linebreak[0] \ 9
\linebreak[0] \ 11 \linebreak[0] \ 7 \linebreak[0] \ 6
\linebreak[0] \ 8 \linebreak[0] \ {\bf 12} \linebreak[0] \ 13
\linebreak[0] \ 5 \linebreak[0] \ 3 \linebreak[0] \ 4
\linebreak[0] \ 2 \linebreak[0] \ 18 \linebreak[0] \ 19
\linebreak[0] \ 1 \linebreak[0] \ 21$.

\begin{corollary}
The number of involutions containing $132$ exactly once of length
$n$ having $p$ fixed points with $1 \leq p \leq n$ is the ballot
number $\binom{n-2}{\frac{n+p}{2}-1}-\binom{n-2}{\frac{n+p}{2}}$.
Moreover, the number of involutions containing $132$ exactly once
of length $n$ is $\binom{n-2}{\lfloor \frac{n-3}{2} \rfloor}$.
\end{corollary}

\begin{proof}
We immediately deduce this result from bijection $\Psi$ of
Theorem~\ref{bij132once-132} and Corollary~\ref{nb132ptfix}.
In fact, the number of involutions containing $132$ once of length $n$
is either the number of $132$-avoiding involutions of length $n-2$
if $n$ is odd or the number of $132$-avoiding involutions of length $n-2$
having more than one fixed point if $n$ is even.
\end{proof}

Of course, some of the following results can immediately be obtained
from bijection $\Psi$ of Theorem~\ref{bij132once-132} and results
of Section~\ref{sec2}.

Let $J_\tau(n)$ denote the number of involutions in $S_n(\tau)$
such containing $132$ exactly once, and let $J_\tau(x)=\sum_{n\geq
0} J_\tau(n)x^n$ be the corresponding generating function. The
following proposition is the base of all the other results in
this section, which holds immediately from definitions.

\begin{proposition}\label{prom2}
Let $\pi$ an involution in $S_n$ such that contains $132$ exactly
once, and let $\pi_j=n$. Then holds either
\begin{enumerate}
\item   or $\pi_n=n$;

\item   or $\pi=(\pi',n,\pi'',\pi''',j)$ where $1\leq j\leq n/2$, $\pi'''={\pi'}^{-1}$ and
    $\pi'$ avoids $132$.

\item   or $\pi=(\pi',m,2m+1,\pi'',m+1)$ where $n=2m+1$, $\pi''={\pi'}^{-1}$
    and $\pi'\in S_{m-1}(132)$.
\end{enumerate}
\end{proposition}

Another approach to find the generating function of involutions in
$S_n$ containing $132$ exactly once is by use Proposition
\ref{prom2}.

\begin{theorem}\label{thcc1}
Let $C(t)$ be the generating function for the Catalan numbers;
then
    $$J_\varnothing(x)=\frac{x^3C(x^2)}{1-x-x^2C(x^2)}.$$
\end{theorem}
\begin{proof}
According to Proposition \ref{prom2} with terms of generating
functions we get the following: the first part of the proposition
yields $xJ_\varnothing(x)$, the second part of the proposition
yields $x^2C(x^2)J_\varnothing(x)$, and the third part of the
proposition gives $x^3C(x^2)$. Hence
    $$J_\varnothing(x)=xJ_\varnothing(x)+x^2C(x^2)J_\varnothing(x)+x^3C(x^2).$$
\end{proof}

\begin{example}\label{excc1}
From Proposition \ref{prom2} it is easy to see that, the number of
the involutions in $S_n(123)$ and containing $132$ exactly once is
$2^{(n-3)/2}$ for $n$ odd, otherwise is $0$. Also,
$J_{1234}(n)=F_{n-3}$ the $(n-3)$th Fibonacci number,
$J_{12345}(n)=3^{[(n-3)/2]}$.
\end{example}

Again, the case of varying $k$ is more interesting. As an
extension of Example \ref{excc1} let us consider the case
$\tau=12\dots k$.

\begin{theorem}\label{thcc2}
For all $k\geq 1$,
    $$J_{12\dots k}(x)=\frac{x}{U_k\left( \frac{1}{2x} \right)} \sum_{j=1}^{k-2} U_{j}\left( \frac{1}{2x} \right).$$
\end{theorem}
\begin{proof}
Proposition \ref{prom2} with use the generating function of permutations in $S_n(132,12\dots k)$
given by $R_k(x)$, yields
    $$J_{12\dots k}(x)=xJ_{12\dots(k-1)}+x^2R_{k-1}(x^2)J_{12\dots k}(x)+x^3R_{k-1}(x^2).$$
By use the relation $R_k(y)=1/(1-yR_{k-1}(y))$ we get that
    $$J_{12\dots k}(x)=xR_k(x^2)J_{12\dots(k-1)}(x)+x^3R_{k-1}(x^2)R_k(x^2),$$
so induction on $k$ with Example \ref{excc1} gives the theorem.
\end{proof}

Similarly,
we obtain another case $\tau=213\dots k$.

\begin{theorem}\label{thcc22}
For all $k\geq 3$,
    $$J_{213\dots k}(x)=\frac{x}{U_k\left( \frac{1}{2x} \right)}
    \left[
    xU_2\left(\frac{1}{2x}\right) +\sum_{j=2}^{k-2} U_{j}\left( \frac{1}{2x} \right)
    \right].$$
\end{theorem}
\begin{proof}
Similarly as proof of Theorem \ref{thcc2} with use the generating
function for the number of $213\dots k$-avoiding permutations in
$S_n(132)$ is given by $R_k(x)$ (see \cite{MV2}), we obtain that
    $$J_{213\dots k}(x)=xR_k(x^2)J_{213\dots(k-1)}(x)+x^3R_{k-1}(x^2)R_k(x^2),$$
and by induction with $J_{213}(x)=x^4R_3(x^2)$ (it is easy to see)
the theorem holds.
\end{proof}

\begin{example}
Theorem \ref{thcc22} yields
$J_{2134}(2n+3)=J_{2134}(2n+4)=F_{2n}$ the $(2n)$th Fibonacci number
for all $n \geq 0$.
\end{example}

\begin{example}\label{excc2}
Proposition \ref{prom2} yields, $J_{231}(n)=1$ for all $n\geq 1$, and
$J_{2341}(n)=2^{[(n-1)/2]}-1$ for all $n\geq 1$.
\end{example}

Once again, the case of varying $k$ is more interesting. As an
extension of Example \ref{excc2} let us consider the case
$\tau=23\dots k1$.

\begin{theorem}\label{thcc3}
For all $k\geq 3$,
    $$J_{23\dots k1}(x)=\frac{x^2U_{k-3}\left( \frac{1}{2x} \right)}{(1-x)U_{k-2}\left( \frac{1}{2x} \right)}
       \left[ 1+\frac{1}{U_{k-1}\left( \frac{1}{2x} \right)}\sum_{j=1}^{k-3} U_j\left( \frac{1}{2x} \right) \right].$$
\end{theorem}
\begin{proof}
Similarly as proof Theorem \ref{thcc2} we have that
    $$J_{23\dots k1}(x)=xJ_{23\dots k1}(x)+x^2R_{k-2}(x^2)J_{12\dots(k-1)}(x)+x^3R_{k-2}(x^2),$$
so by using Theorem \ref{thcc2} the theorem holds.
\end{proof}

More generally, we present an explicit expression when
$\tau=[k,d]$ as follows.

\begin{theorem}\label{thcc4}
For $k\geq 4$, $2\leq d\leq k/2$,
$$J_{[k,d]}(x)=\frac{R_d(x^2)}{1-xR_d(x^2)}\left[
    x^2R_{k-d-1}(x^2)+\frac{x^2(R_{k-d-1}(x^2)-R_{d-1}(x^2))}{U_{k-d}\left(\frac{1}{2x}\right)}\sum_{j=1}^{k-d-2} U_j\left(\frac{1}{2x}\right)
    \right].$$
\end{theorem}
\begin{proof}
According to Proposition \ref{prom2} in terms of generating
functions we get the following. In first case $xJ_{[k,d]}(x)$. In
the third case, if $\pi'$ contains $12\dots(k-d-1)$ then $\pi$
contains $[k,d]$ which is a contradiction, we get that $\pi'$
avoids $12\dots(k-d-1)$, hence $x^3R_{k-d-1}(x^2)$. Finally, in
the second case, let us observe two subcases $\pi''$ contains
$12\dots (k-d)$ or avoids $12\dots (k-d)$; so by use the same
argument of the third case we get
$$x^2R_{k-d-1}(x)J_{12\dots (k-d)}(x)+x^2R_{d-1}(x^2)(J_{[k,d]}(x)-J_{12\dots(k-d)}(x)).$$
Therefore, if we add all these cases we get $J_{[k,d]}(x)$. Hence,
by Theorem \ref{thcc2} this theorem holds.
\end{proof}

%
\section{Containing $132$ once and containing another pattern}
Let $J_\tau^r(n)$ denote the number of involutions in $S_n$ such
containing $132$ exactly once and containing $\tau$ exactly $r$
times. Let $J_\tau^r(x)=\sum_{n\geq 0} J_\tau^r(n)x^n$ be the
corresponding generating function. Let us start be the following
result.

\begin{theorem}\label{thd1}
For all $k\geq 1$,
   $$J_{12\dots k}^1(x)=0.$$
\end{theorem}
\begin{proof}
By Proposition \ref{prom2} it is easy to see
    $$J_{12\dots k}^1(x)=xJ_{12\dots (k-1)}^1(x)+x^2R_{k-1}(x^2)J_{12\dots k}(x).$$
with $J_{12}^1(x)=0$, hence induction on $k$ gives the theorem.
\end{proof}

Similarly as Theorem \ref{thd1} we have another case where
$\tau=23\dots k1$.

\begin{theorem}\label{thd2}
For all $k\geq 1$,
    $$J_{23\dots k1}^1(x)=0.$$
\end{theorem}

\begin{example}\label{exdd1}
Proposition \ref{prom2} yields the following. The number of
involutions $J_{21}^1(n)=1$ for $n\geq 3$, and
$J_{213}(n)=2^{(n-8)/2}(1+(-1)^n)$.
\end{example}

Once again, the case of varying $k$ is more interesting. As an
extension of Example \ref{exdd1} let us consider the case
$\tau=213\dots k$.

\begin{theorem}\label{thd3}
For all $k\geq 3$, $J_{213\dots
    k}^1(x)=\frac{x(1-x^2)}{U_k\left(\frac{1}{2x}\right)}$.
\end{theorem}
\begin{proof}
By Proposition \ref{prom2} it is easy to obtain
    $$J_{213\dots k}^1(x)=xJ_{213\dots (k-1)}^1(x)+x^2R_{k-1}(x^2)J_{213\dots k}(x).$$
with $J_{21}^1(x)=x^3/(1-x)$ (which is yield directly from
definitions), hence induction on $k$ gives the theorem.
\end{proof}

%


\end{document}